\documentclass{article}

\usepackage{latexsym,amssymb,graphicx,amscd}
\usepackage{amsfonts}

\usepackage{amsthm}
\usepackage{epsf}
\usepackage{color}

\newtheorem{theorem}{Theorem}[section]    
\newtheorem{lemma}[theorem]{Lemma}

\newtheorem{proposition}[theorem]{Proposition}

\def\PSL{\mbox{\rm{PSL}}}
\def\P{\mbox{\rm{P}}}

\def\SL{\mbox{\rm{SL}}}
\def\SO{\mbox{\rm{SO}}}
\def\SU{\mbox{\rm{SU}}}
\def\Sp{\mbox{\rm{Sp}}}
\def\PSU{\mbox{\rm{PSU}}}

\def\GL{\mbox{\rm{GL}}}

\def\Out{\mbox{\rm{Out}}}

\def\qed{ $\sqcup\!\!\!\!\sqcap$}

\usepackage{amscd}

\newcommand{\BZ}{{\bf{Z}}}

\newcommand{\BF}{{\bf{F}}}

\title{Frattini and related subgroups of Mapping Class Groups}
\author{G. Masbaum 
~and A. W. Reid \thanks{This work was partially supported by the
  N. S. F.}}

\date{January 14, 2015}

\begin{document} 

\maketitle

\centerline{To V. P. Platonov on the occasion of his 75th birthday}

\begin{abstract} Let $\Gamma_{g,b}$ denote the orientation-preserving
Mapping Class Group of a closed orientable surface of genus $g$ with
$b$ punctures.  For a group $G$ let $\Phi_f(G)$ denote the
intersection of all maximal subgroups of finite index in $G$.
Motivated by a question of Ivanov as to whether $\Phi_f(G)$ is
nilpotent  
when $G$ is a finitely generated subgroup 
of $\Gamma_{g,b}$,
in this paper we 
compute
$\Phi_f(G)$ for certain
subgroups of $\Gamma_{g,b}$. In particular, we answer Ivanov's question
in the affirmative
for 
these
subgroups of $\Gamma_{g,b}$.
\end{abstract}
          

\section{Introduction}
We fix the following notation throughout this paper:   Let $\Gamma_{g,b}$
denote the orientation-preserving Mapping Class Group of a closed orientable surface of genus $g$ with $b$ punctures.  When $b=0$ we simply write
$\Gamma_g$.  In addition, when $b>0$
we let $\P\Gamma_{g,b}$ denote the {\em pure Mapping Class Group};
i.e. the subgroup of $\Gamma_{g,b}$ consisting of those
elements that fix the punctures pointwise. 
The Torelli group ${\cal I}_g$ is the subgroup of $\Gamma_g$ arising as
the kernel of the
homomorphism $\Gamma_g \rightarrow \Sp(2g,{\bf Z})$ 
coming from the action of $\Gamma_g$ on $H_1(\Sigma_g,{\bf Z})$. As usual  $\Out(F_n)$ will denote the Outer Automorphism Group of a free group
of rank $n\geq 2$. 

For a group $G$, the {\em Frattini
  subgroup} $\Phi(G)$ of $G$ is defined to be the intersection of all
maximal 
subgroups of $G$ (if they exist), otherwise it is defined to
be the group $G$ itself. 
(Here a maximal subgroup is a strict subgroup which is maximal with respect to inclusion.)  
 In addition we define $\Phi_f(G)$ to be the
intersection of all maximal subgroups of finite index in $G$. Note that $\Phi(G) < \Phi_f(G)$.\footnote{We use the notation $G_1<G_2$ to indicate that $G_1$ is a subgroup of $G_2$ (including the case where $G_1=G_2$).}

Frattini's original theorem is that if $G$ is finite then $\Phi(G)=\Phi_f(G)$ is a nilpotent group (see
for example \cite[Theorem 11.5]{Ro}). For infinite groups this is not the case: there are examples of finitely generated infinite groups $G$ with $\Phi(G)$ not nilpotent \cite[p.~328]{H}. 
On the other hand,  
in \cite{Pl}
Platonov showed that if $G$ is any finitely generated linear group then
$\Phi(G)$ and $\Phi_f(G)$ are nilpotent.

Motivated by the question as to whether $\Gamma_g$ is a linear group,
in \cite{Lo}, Long proved that $\Phi(\Gamma_g)=1$ for $g\geq 3$, and
$\Phi(\Gamma_2)={\bf Z}/2{\bf Z}$.  This was extended by Ivanov
\cite{Iv1} who showed that (as in the linear case), $\Phi(G)$ is
nilpotent for any finitely generated subgroup
$G<\Gamma_{g,b}$. 
Regarding $\Phi_f(G)$, in \cite{Iv1} (and then again
in \cite{Iv2}), Ivanov asks whether the same is true for $\Phi_f$: 

\vskip 8pt

\noindent{\bf Question:} (Ivanov \cite{Iv1,Iv2}) Is $\Phi_f(G)$ is nilpotent for every
finitely generated subgroup $G$ of $\Gamma_{g,b}$? \\[\baselineskip]
The aim of this note is to prove some results in the
direction of answering Ivanov's question. 
In particular, the following theorem answers Ivanov's question in the affirmative for $\Gamma_g$ and some of its subgroups in the case where $g\geq 3$.

\begin{theorem}
\label{main1}
Suppose that $g\geq 3$, and that $G$ is either 
\begin{enumerate}
\item[(i)] the Mapping Class Group $\Gamma_g$, or
\item[(ii)] a normal subgroup of $\Gamma_g$ (for example the Torelli group ${\cal I}_g$, the Johnson kernel ${\cal K}_g$, or any higher term in the Johnson filtration of $\Gamma_g$), or 
\item[(iii)] a subgroup of $\Gamma_g$ which contains a finite index subgroup of  the Torelli group  ${\cal I}_g$.
\end{enumerate} 
Then $\Phi_f(G)=1$.\end{theorem}

\vskip 8pt

\noindent{\bf Remarks:} (1)~Since $\Phi(G) < \Phi_f(G)$, 
our methods also give a different proof of Long's result that 
$\Phi(\Gamma_g)=1$ for $g\geq 3$.  As for the case $g\leq 2$, note
that $\Gamma_1$ and $\Gamma_2$ are linear (see \cite{BB} for $g=2$) and so 
Platonov's result \cite{Pl} applies to 
answer Ivanov's question in the affirmative 
in these cases for all finitely
generated subgroups. 
On the other hand, for $g\geq 3$ the Mapping Class Group is not known to be linear and no other technique for answering Ivanov's question was known. In fact,  as pointed out by Ivanov, neither  the methods of
\cite{Lo} or \cite{Iv1} apply to  $\Phi_f$, and so even the case of
$\Phi_f$ of the Mapping Class Group itself was not known. \\[\baselineskip] 
\noindent (2)~Note that
in Platonov's and Ivanov's theorems and in Ivanov's question, 
the Frattini subgroup and its variant $\Phi_f$ are 
considered 
for finitely generated subgroups.  In reference to Theorem
\ref{main1} above, it remains an open question as to whether the
Johnson kernel ${\cal K}_g$, or any higher term in the Johnson
filtration of $\Gamma_g$, is finitely generated or not.\\[\baselineskip]
Perhaps the most interesting feature of the proof of Theorem \ref{main1} is 
that it is another application of the
projective unitary representations arising 
in Topological Quantum Field Theory (TQFT) first constructed by Reshetikhin and
Turaev \cite{RT} (although as in \cite{MR}, the perspective here is that of
the skein-theoretical approach of \cite{BHMV}).

We are also able to prove:

\begin{theorem}
\label{main1add}
Assume that $b>0$, then
$\Phi_f(\P\Gamma_{g,b})$ is either trivial or ${\bf Z}/2{\bf Z}$. Indeed,
$\Phi_f(\P\Gamma_{g,b})=1$ unless $(g,b)\in\{(1,b),(2,b)\}$. 
\end{theorem}

\noindent The reason for separating out the case when $b>0$ is that
the proof does not directly use the TQFT framework, but rather makes use of
Theorem \ref{main1}(i) in conjunction with the Birman exact sequence
and a general group theoretic lemma (see 
      \S \ref{sec6}).
  We expect that the methods
of this paper will also answer Ivanov's question for $\Gamma_{g,b}$ but at
present we are unable to do so. We comment further on this 
at the end of  \S \ref{sec6}.

In addition our methods can also be used to give a straightforward proof of the following.

\begin{theorem}
\label{main2}
Suppose that $n\geq 3$, $\Phi(\Out(F_n))=\Phi_f(\Out(F_n))=1$.\end{theorem}

\noindent Note that it was shown in \cite{Hu} that $\Phi(\Out(F_n))$ is 
finite.

As remarked upon above, this note was largely motivated by the
questions of Ivanov.  To that end, we discuss a possible approach to
answering Ivanov's question in general using the aforementioned
projective unitary representations arising from TQFT, coupled with 
Platonov's work \cite{Pl}.  Another motivation for this work arose
from attempts to understand the nature of the Frattini subgroup and
the center of the profinite completion of $\Gamma_g$
and ${\cal I}_g$.  We discuss  these further in 
      \S \ref{sec7}.\\[\baselineskip]
\noindent{\bf Acknowledgements:}~{\em The authors wish to thank the organizers of the conference "Braids and Arithmetic" at CIRM Luminy in October 2014, where this work was completed.}

\section{Proving triviality of $\Phi_f$}

Before stating and proving an elementary but
useful technical result we introduce some notation.
Let $\Gamma$ be a finitely generated group, and let ${\cal S}=\{G_n\}$
a collection of finite groups together with epimorphisms 
$\phi_n:\Gamma\rightarrow G_n$.
We say that $\Gamma$ is 
{\em residually}-${\cal S}$, if given any non-trivial element
$\gamma\in\Gamma$, there is some group $G_n\in{\cal S}$
and an epimorphism $\phi_n$ for which $\phi_n(\gamma)\neq 1$.
Note that, as usual, this is equivalent to the statement 
$\bigcap \ker\phi_n = 1$.

\begin{proposition}
\label{tool}
Let $\Gamma$ and ${\cal S}$ be as above with $\Gamma$ being
residually-${\cal S}$. Assume further that $\Phi(G_n)=1$ for every 
$G_n\in{\cal S}$. Then $\Phi(\Gamma)=\Phi_f(\Gamma)=1$.\end{proposition}

Before commencing with the proof of this proposition, we recall the following
property.

\begin{lemma}
\label{frattini_under_epi}
Let $\Gamma$ and $G$ be groups and $\alpha:\Gamma\rightarrow G$ an epimorphism.
Then $\alpha(\Phi(\Gamma)) \subset \Phi(G)$ and 
$\alpha(\Phi_f(\Gamma)) \subset \Phi_f(G)$.\end{lemma}

\noindent{\bf Proof:}~We prove the last statement. Let $M$ be a maximal
subgroup of $G$ of finite index. Then $\alpha^{-1}(M)$ is a maximal subgroup
of $\Gamma$ of finite index in $\Gamma$, and hence $\Phi_f(\Gamma) \subset 
\alpha^{-1}(M)$. Thus $\alpha(\Phi_f(\Gamma)) \subset M$ for all maximal
subgroups $M$ of finite index in $G$ and the result 
follows.\qed\\[\baselineskip]
\noindent{\bf Remark:}~As pointed out in \S 1,  for a finite group $G$, 
$\Phi(G)=\Phi_f(G)$.\\[\baselineskip]
\noindent{\bf Proof of Proposition \ref{tool}:}~We give the argument for
$\Phi_f(\Gamma)$, the argument for $\Phi(\Gamma)$ is exactly the same.
Thus suppose that $g\in\Phi_f(\Gamma)$ is a non-trivial element. Since
$\Gamma$ is residually-${\cal S}$, there exists some $n$ so that
$\phi_n(g) \in G_n$ is non-trivial.  
However, by Lemma \ref{frattini_under_epi} (and the remark following it) we 
have:

$$\phi_n(g) \in \phi_n(\Phi_f(\Gamma)) < \Phi_f(G_n)=\Phi(G_n),$$ 
and in particular $\Phi(G_n)\neq 1$, a contradiction.\qed

\section{The quantum representations and finite quotients}

We briefly recall some of \cite{MR} (which uses \cite{BHMV} and \cite{GM}).
As in \cite{MR} we only consider the case of $p$ a prime satisfying
$p\equiv 3 \pmod 4 $. 

Let $\Sigma$ be  a
closed orientable surface of genus $g\geq 3$.
The integral $SO(3)$-TQFT constructed in
\cite{GM} provides a representation of a central extension $\widetilde
\Gamma_g$  of $\Gamma_g$ 
$$\rho_p \,:\, \widetilde
\Gamma_g \longrightarrow \GL(N_g(p),\BZ[\zeta_p])~,$$ where
$\zeta_p$ is a primitive $p$-th root of unity, 
$\BZ[\zeta_p]$ 
is
the ring of cyclotomic integers and 
$N_g(p)$ the dimension of a vector space 
$V_p(\Sigma)$
on which the representation
acts. It is known that $N_g(p)$ 
is given by a Verlinde-type formula and goes to infinity as
$p\rightarrow \infty$. For convenience we simply set $N=N_g(p)$.

As in \cite{MR} the
image group $\rho_p(\widetilde{\Gamma}_g)$ will be denoted by
$\Delta_g$.  As is pointed out in \cite{MR}, 
$\Delta_g< \SL(N, \BZ[\zeta_p])$, and moreover, 
$\Delta_g$ is actually contained in a special unitary group
$\SU(V_p,H_p;\BZ[\zeta_p])$, where $H_p$ is a Hermitian form defined over 
the real field 
${\bf Q}(\zeta_p+\zeta_p^{-1})$.

Furthermore, the homomorphism $\rho_p$, descends to a projective
representation of $\Gamma_g$ (which we
denote by $\overline{\rho}_p$):

$$\overline{\rho}_p : \Gamma_g \longrightarrow \PSU(V_p,H_p;\BZ[\zeta_p]),$$

What we need from \cite{MR} is the following.
We can find infinitely many rational primes $q$ which split completely in
$\BZ[\zeta_p]$, and for every such prime $\tilde q$ of $\BZ[\zeta_p]$ lying
over such a $q$, we can consider the group 
$$\pi_{\tilde q}(\Delta_g) \subset \SL(N,q),$$ where $\pi_{\tilde q}$ is the reduction homomorphism from
$\SL(N,\BZ[\zeta_p])$ to $\SL(N,q)$ induced by the isomorphism
$\BZ[\zeta_p]/\tilde q\simeq \BF_q$.
As is shown in \cite{MR} (see also \cite{Fu})
we obtain epimorphisms $\Delta_g\twoheadrightarrow \SL(N,q)$ for
all but finitely many of these primes $\tilde{q}$, and it then 
follows easily that we obtain epimorphisms
$\Gamma_g \twoheadrightarrow \PSL(N,q)$.  We denote these epimorphisms by
$\rho_{p,\tilde{q}}$.  These should be thought of as reducing 
the images of 
      $\overline{\rho}_{p}$
modulo $\tilde{q}$. That one obtains finite simple groups of the form $\PSL$
rather than $\PSU$ when $q$ is a split prime is discussed in \cite{MR} \S2.2.

\begin{lemma}
\label{residual}
For each $g\geq 3$, $\bigcap \ker\rho_{p,\tilde{q}}=1$.\end{lemma}

\noindent{\bf Proof:}~Fix $g\geq 3$ and suppose that there exists a
non-trivial element $\gamma\in \bigcap \ker\rho_{p,\tilde{q}}$. Now it
follows from asymptotic faithfulness \cite{A,FWW}
that $\bigcap \ker\overline{\rho}_p = 1$. Thus
for some $p$ there exists $\overline{\rho}_p$ such that 
$\overline{\rho}_p(\gamma)\neq
1$. Now $\rho_{p,\tilde{q}}(\gamma)$ is obtained by reducing
$\overline{\rho}_p(\gamma)$ modulo $\tilde{q}$, and so there clearly exists
$\tilde{q}$ so that
$\rho_{p,\tilde{q}}(\gamma)\neq 1$, a contradiction.\qed\\[\baselineskip]

\section{Proofs of Theorems \ref{main1} and \ref{main2}}

The proof of Theorem \ref{main1} 
for $G=\Gamma_g$ 
follows easily as a special case of our next result. To state this, we introduce some notation:  If $H<\Gamma_g$, we denote by $\widetilde H$, the inverse image of $H$ under the projection $\widetilde{\Gamma}_g \rightarrow \Gamma_g$.

\begin{proposition}
\label{saturating_trivial_frattini}
Let $g\geq 3$, and assume that $H$ is a finitely generated subgroup
of $\Gamma_g$ for which $\rho_p(\widetilde H)$ has the same Zariski closure and adjoint trace field as  $\Delta_g$. Then $\Phi(H)=\Phi_f(H)=1$.\end{proposition}

\noindent{\bf Proof:}~We begin with a remark.  That the homomorphisms $\rho_{p,\tilde{q}}$ of \S 3 are
surjective is proved using Strong Approximation. The main ingredients of this are the Zariski density of $\Delta_g$ in the algebraic group $\SU(V_p,H_p)$, and the fact that the adjoint trace field of $\Delta_g$ is the field ${\bf Q}(\zeta_p+\zeta_p^{-1})$ over which the group  $\SU(V_p,H_p)$ is defined (see \cite{MR} for more details). In particular, the proof establishes surjectivity  of $\rho_{p,\tilde{q}}$ when restricted to any subgroup $H<\Gamma_g$  equipped with the hypothesis of the proposition.

To complete the proof,  the groups $\PSL(N,q)$ are
finite simple groups (since the dimensions $N$ are all very large)
so their Frattini subgroup is trivial. 
This follows from Frattini's theorem, or, more simply, from the fact that the Frattini subgroup of a finite group is a normal subgroup which is moreover a strict subgroup (since finite groups do have maximal subgroups). 
Hence the result follows from
Lemma \ref{residual},
Proposition \ref{tool} and the remark at the start of the proof.\qed\\[\baselineskip]
\noindent In particular,   $\Gamma_g$ satisfies the hypothesis of Proposition \ref{saturating_trivial_frattini},  and so  $\Phi_f(\Gamma_g)=1$. This
also recovers the result of Long \cite{Lo} proving triviality of the
Frattini subgroup.\\[\baselineskip]
\noindent The proof of Theorem \ref{main1} in case (ii), that is, when $G$ is a normal subgroup of $\Gamma_g$, follows from this and  the following general fact:

\begin{proposition}\label{fFN} If $N $ is a normal subgroup of a group $\Gamma$, then  $\Phi_f(N) < \Phi_f(\Gamma)$. 
\end{proposition}

This fact is known for Frattini subgroups of finite groups, and the proof can be adapted to our situation. We defer the details to Section~\ref{FFN}.\\[\baselineskip]
In the remaining case (iii) of Theorem \ref{main1}, $G$  is a subgroup of $\Gamma_g$ which contains a finite index subgroup of  the Torelli group ${\cal I}_g$. We shall show that $G$ satisfies the
hypothesis of  Proposition \ref{saturating_trivial_frattini}, and deduce $\Phi_f(G)=1$ as  before. 

Consider first the case where $G$ is the Torelli group ${\cal I}_g$ itself. 
Recall the short exact sequence 

$$1 \longrightarrow {\cal I}_g \longrightarrow \Gamma_g \longrightarrow \Sp(2g,{\bf Z})\longrightarrow 1~.$$
 We now use the following well-known facts. 

\begin{enumerate}
\item[$\bullet$] $\Gamma_g$ is generated by Dehn twists, which map to transvections in $ \Sp(2g,{\bf Z})$.
\item[$\bullet$] The central extension $\widetilde \Gamma_g$ of $\Gamma_g$ is generated by certain lifts of Dehn twists, and $\rho_p$ of every such lift is a matrix of order $p$.
\item[$\bullet$]  The quotient of $\Sp(2g,{\bf Z})$ by the normal subgroup generated by $p$-th powers of transvections is the finite group $\Sp(2g,{\bf Z}/p{\bf Z})$ (see \cite{BMS} for example). 
\end{enumerate}

\noindent It follows that the finite group $\Sp(2g,{\bf Z}/p{\bf Z})$ 
admits a surjection 
onto
the quotient group  
$$\Delta_g / \rho_p(\widetilde{\cal I}_g)$$
(recall that $\Delta_g= \rho_p(\widetilde{\Gamma}_g)$) and hence the group $\rho_p(\widetilde{\cal I}_g)$ has finite index in $\Delta_g$. 
But the Zariski closure of $\Delta_g$ is the connected, simple, algebraic group $\SU(V_p,H_p)$.
Thus $\rho_p(\widetilde{\cal I}_g)$ and  $\Delta_g$ have the same Zariski closure.  Again using the fact that $\SU(V_p,H_p)$ is a simple algebraic group,
we also deduce that $\rho_p(\widetilde{\cal I}_g)$ has the same adjoint trace field as $\Delta_g$ (this follows from  \cite{DM} Proposition 12.2.1 for example).
This shows that ${\cal I}_g$ indeed satisfies the hypothesis of Proposition~\ref{saturating_trivial_frattini}, and so once again
$\Phi_f({\cal I}_g)=1$. The same arguments work when $G$ has finite index in  ${\cal I}_g$, and also when $G$ is any subgroup of $\Gamma_g$ which contains a finite index subgroup of ${\cal I}_g$. This completes the proof of Theorem \ref{main1}.\qed\\[\baselineskip]
We now turn to the proof of Theorem \ref{main2}. 
To deal with the case of $\Out(F_n)$, we recall that R.~Gilman \cite{Gil}
showed that for $n\geq 3$, $\Out(F_n)$
is residually alternating: i.e. in the notation of \S 2, the 
collection ${\cal S}$ consists of alternating groups.\\[\baselineskip]
\noindent{\bf Proof of Theorem \ref{main2}:}~For $n\geq 3$, the abelianization
of $\Out(F_n)$ is ${\bf Z}/2{\bf Z}$ (as can be seen directly from Nielsen's 
presentation of $\Out(F_n)$, see \cite{Vo} \S 2.1). Hence,
$\Out(F_n)$ does not admit a surjection 
onto
$A_3$ or $A_4$. Thus all the alternating quotients described by Gilman's
result above have trivial Frattini subgroups (as in the proof of 
Proposition \ref{saturating_trivial_frattini}).  The proof is completed
using the residual alternating property and Proposition \ref{tool}.\qed

\section{Proof of Proposition~\protect{\ref{fFN}}}\label{FFN}

Let  $\Gamma$ be a group and $N $ a normal subgroup of $\Gamma$. We wish to show that $\Phi_f(N) < \Phi_f(\Gamma)$. We proceed as follows.

 First a preliminary observation. Let $K=\Phi_f(N)$. It is easy to see that $K$ is characteristic in $N$ (i.e., fixed by every automorphism of $N$).  Since $N$ is normal in $\Gamma$, it follows that  $K$ is normal in $\Gamma$. This implies that for every subgroup $M$ of $\Gamma$, the set $$KM=\{km\,|\, k\in K, m\in M\}$$ is a subgroup of $\Gamma$. Moreover, since $K < N$, we have 

\begin{equation}\label{eq1} KM \cap N = KM_1\end{equation} 
where $M_1= M\cap N$. 
To see the inclusion $KM \cap N \subset KM_1$, write an element of $KM \cap N$ as $km=n$ and observe that $m\in N$ since  $K<N$. Thus $m\in M_1$. The 
reverse inclusion is immediate.

Now suppose  for a contradiction that $K =\Phi_f(N)$ is not contained in $ \Phi_f(\Gamma)$. Then there 
exists
a maximal subgroup $M<\Gamma$ of finite index such that $K$ is not contained in $M$. Write $$M_1=M\cap N$$ as above. Then $M_1$ is a finite index subgroup of $N$. If $M_1=N$ then $N$ is contained in $M$, and hence so is $K$, which is a contradiction.  Thus $M_1$ is a strict subgroup of $N$, and since its index in $N$ is finite, $M_1$ is contained in a maximal subgroup $H$, say, of $N$. 

The proof is now concluded as follows. By definition, $K=\Phi_f(N)$ is also contained in $H$. Hence the group $KM_1$ is contained in $H$ and therefore strictly smaller than $N$. On the other hand, by the maximality of $M$ in $\Gamma$, we have $KM=\Gamma$,  and hence, using (\ref{eq1}), we have 
$$ KM_1 = KM \cap N =\Gamma \cap N = N~.$$ This contradiction completes the proof.\qed\\[\baselineskip]
\noindent{\bf Remark:} If we consider the original Frattini group $\Phi$ in place of $\Phi_f$, one can show similarly that $\Phi(N) < \Phi(G)$, provided that every subgroup of $N$ is contained
in a maximal subgroup of $N$; e.g. when $N$ is finitely generated.

\section{Proof of Theorem \ref{main1add}}\label{sec6}

We begin by recalling the Birman exact sequence. Let $\Sigma_{g,b}$ denote the closed orientable surface of genus $g$ with $b$ punctures. If $b=0$ we abbreviate to $\Sigma_g$. There is a short exact sequence (the
{\em Birman exact sequence}):

$$1\rightarrow \pi_1(\Sigma_{g,(b-1)})\rightarrow \P\Gamma_{g,b}\rightarrow \P\Gamma_{g,(b-1)}\rightarrow 1,$$

\noindent where the map $\P\Gamma_{g,b}\rightarrow \P\Gamma_{g,(b-1)}$ is the
forgetful map, and the map $\pi_1(\Sigma_{g,(b-1)})\rightarrow \P\Gamma_{g,b}$ the
point pushing map (see \cite{FM} Chapter 4.2 for details). Also
in the case when $b=1$, the symbol $\P\Gamma_{g,0}$ simply denotes the
Mapping Class Group $\Gamma_g$.

It will be useful to recall that an alternative description of
$\P\Gamma_{g,b}$ is as the kernel of an epimorphism 
$\Gamma_{g,b}\rightarrow S_b$ (the symmetric group on $b$ letters).

The proof will proceed by induction, using Theorem \ref{main1}(i) to
get started, together with the following (which is an adaptation of 
Lemma 3.5 of \cite{ABetal} to the case of $\Phi_f$). The
proof is included 
in \S \ref{sec.all}
below. 

We introduce the following
notation. Recalling \S 2, let $G$ be a group, say $G$ is {\em
  residually simple} if the collection ${\cal S}=\{G_n\}$ (as in \S 2)
consists of finite non-abelian simple groups.

\begin{lemma}
\label{allenby}
Let $N$ be a finitely generated normal subgroup of the group $G$ and assume that $N$ is residually simple. Then $N\cap \Phi_f(G)=1$. In particular if $\Phi_f(G/N)=1$, then $\Phi_f(G)=1$.\end{lemma}

Given this we now complete the proof. In the cases of 
$(0,1)$, $(0,2)$ and $(0,3)$, it is easily seen that
the subgroup $\Phi_f$ is trivial. 

Thus we now assume that we are not in those cases. As is
well-known, $\pi_1(\Sigma_{g,b})$ is residually simple for those
surface groups under consideration, except the case of $\pi_1(\Sigma_1)$
which we deal with separately below. For example this follows by
uniformization of the surface by a Fuchsian group with algebraic
matrix entries and then use Strong Approximation.

Assume first that $g\geq 3$, then 
Theorem \ref{main1}(i), Lemma \ref{allenby} and the Birman
exact sequence immediately proves that
the statement holds for $\P\Gamma_{g,1}$. The remarks above, Lemma \ref{allenby}
and induction then proves the result for $\P\Gamma_{g,b}$ whenever $g\geq 3$
and $b>0$.

Now assume that $g=0$. The base case of the induction here is
$\P\Gamma_{0,4}$. From 
the
above, it is easy to see that $\P\Gamma_{0,3}$ is trivial,
and so $\P\Gamma_{0,4}$ is a free group of rank $2$. As such, it follows that
$\Phi_f(\P\Gamma_{0,4})=1$. The remarks above, Lemma \ref{allenby}
and induction then proves the result for $\P\Gamma_{0,b}$ whenever 
$b>0$.

When $g=1$, $\Gamma_1\cong \Gamma_{1,1}\cong \SL(2,{\bf Z})$ and it is easy to
check that $\Phi_f(\SL(2,{\bf Z}))={\bf Z}/2{\bf Z}$ (coinciding
with the center of $\SL(2,{\bf Z})$). Now $\P\Gamma_{1,1}=\Gamma_{1,1}$
and so these facts together with
Lemma \ref{allenby}
and induction then prove the result (i.e. that $\Phi_f(\P\Gamma_{1,b})$
is either trivial or ${\bf Z}/2{\bf Z}$).

In the case of $g=2$, by \cite{BB} $\Gamma_2$ is linear, and so \cite{Pl}
also proves that $\Phi_f(\Gamma_2)$ is nilpotent. We claim that this
forces $\Phi_f(\Gamma_2)={\bf Z}/2{\bf Z}$.  To see this we
argue as follows.

If 
$\Phi_f(\Gamma_2)$ is finite, it is central by \cite{Lo} Lemma 2.2. 
Since $\Phi_f(\Gamma_2)$ contains $\Phi(\Gamma_2)$, 
which is equal to the center ${\bf Z}/2{\bf Z}$ of $\Gamma_2$ by \cite{Lo} Theorem 3.2, it follows 
that $\Phi_f(\Gamma_2) = {\bf Z}/2{\bf Z}$.

Thus it is enough to show that $\Phi_f(\Gamma_2)$ is finite. 
Assume that
it is not. Then by
 \cite{Lo} Lemma 2.5,
$\Phi_f(\Gamma_2)$ contains a pseudo-Anosov element. Indeed, \cite{Lo}  Lemma
2.6 shows that the set of invariant laminations 
of pseudo-Anosov elements
in $\Phi_f(\Gamma_2)$
is dense in projective measured lamination space. This contradicts 
$\Phi_f(\Gamma_2)$ being nilpotent (e.g. the argument of \cite{Lo}
p. 86 constructs a free subgroup).

As before, using Lemma \ref{allenby} and by induction via the Birman
exact sequence,
we can now handle the cases of $\Gamma_{2,b}$ with $b>0$.\qed\\[\baselineskip]
\noindent{\bf Remark 1:}~Recall that
the {\em hyperelliptic
Mapping Class Group} (which we denote by $\Gamma_g^h$)
is defined to be the subgroup of $\Gamma_g$
consisting of those elements that commute with a fixed hyperelliptic 
involution.  It is pointed out in \cite{BB} p. 706, that the
arguments used in \cite{BB} prove that $\Gamma_g^h$ is linear. Hence
once again 
$\Phi_f(G)$ is nilpotent for every finitely generated subgroup $G$ 
of $\Gamma_g^h$.\\[\baselineskip]
\noindent{\bf Remark 2:}~We make some
comments on the case of $\Gamma_{g,b}$ with $b>0$.  First, since
$\P\Gamma_{g,b}=\ker\{\Gamma_{g,b}\rightarrow S_b\}$ and $\Phi_f(S_b)=1$, if
$\P\Gamma_{g,b}$ were known to be residually simple then the argument in the
proof of Theorem \ref{main1add} could be used to show that 
$\Phi_f(\Gamma_{g,b})=1$. Hence we raise here:\\[\baselineskip]
\noindent{\bf Question:}~{\em Is $\P\Gamma_{g,b}$ residually 
simple?}\\[\baselineskip]
Another approach to showing that $\Phi_f(\Gamma_{g,b})=1$ is to directly use the
representations arising from TQFT.  In this case the result of Larsen
and Wang \cite{LW} that allows us to prove Zariski density in \cite{MR}
needs to be established.  Given this, the proof (for most $(g,b)$) would then
follow as above.

\section{Proof of Lemma~\protect{\ref{allenby}}}\label{sec.all}

As already mentioned, in what follows we adapt  the proof of 
Lemma 3.5 of \cite{ABetal} to the case of $\Phi_f$.

We argue by contradiction
and assume that there exists a non-trivial element $x\in N\cap
\Phi_f(G)$. By the residually simple assumption, we can find a non-abelian
finite simple group $S_0$ and an epimorphism $f:N\rightarrow S_0$ for which
$f(x)\neq 1$. Set $K_0=\ker~f$ and let $K_0,K_1\ldots, K_n$ be the
distinct copies of $K_0$ which arise on mapping $K_0$ under the
automorphism group of $N$ (this set being finite since $N$ is finitely
generated). Set $K=\bigcap K_i$, a characteristic subgroup of
finite index in $N$.  As in \cite{ABetal}, it follows from standard
finite group theory that $N/K$ is isomorphic to a direct 
product
of finite
simple groups (all of which are isomorphic to $S_0=N/K_0$).

Now $K$ being characteristic in $N$ implies that $K$ is a normal
subgroup of $G$. Put $G_1=G/K$ and let $f_1:G\rightarrow G_1$ denote the canonical
homomorphism. Also write $N_1$ for $N/K=f_1(N)$.  Now $f_1(x) \in N_1$ and $f_1(x)\in f_1(\Phi_f(G))$
which by Lemma \ref{frattini_under_epi} implies that $f_1(x)\in
\Phi_f(G_1)$. Hence $$f_1(x) \in N_1\cap \Phi_f(G_1)~.$$  Following
\cite{ABetal}, let $C$ denote the centralizer of $N_1$ in $G_1$, and as in
\cite{ABetal}, 
we can deduce various properties about the groups $N_1$
and $C$. Namely:

\medskip

\noindent (i) since $N_1$ is a finite group, its  centralizer $C$ in $G_1$ is of finite index in $G_1$ .

\smallskip

\noindent (ii) since $N_1$ is a product of non-abelian finite simple groups 
it has trivial center, and so $C\cap N_1=1$.

\smallskip

\noindent (iii) since $N_1$ is normal in $G_1$, $C$ is normal in $G_1$.

\medskip

Using (iii), put $G_2=G_1/C$ and let $f_2:G_1\rightarrow G_2$ denote the canonical
homomorphism. Also write $N_2$ for $f_2(N_1)$. Arguing as before (again invoking Lemma \ref{frattini_under_epi}), we have 
$$f_2(f_1(x)) \in N_2\cap \Phi_f(G_2)~.$$ Moreover, since $f_1(x)\in N_1$ and $f_1(x)\neq 1$ by construction, we have  from (ii) that $f_2f_1(x)\neq 1$. Thus the intersection $$H:= N_2\cap \Phi_f(G_2)$$ is a non-trivial group. 

As in \cite{ABetal}, we will now get a contradiction by showing that $H$ is both a nilpotent group and a direct product of non-abelian finite simple groups, which 
is possible only if $H$ is trivial. Here is the argument.  From (i) above we deduce that $G_2$ is a finite group, hence  $\Phi_f(G_2)=\Phi(G_2)$ is  
nilpotent  by Frattini's theorem. Thus $H<\Phi_f(G_2)$ is  nilpotent. On the other hand, $N_2$ is a quotient of $N_1$ and hence a direct product of non-abelian finite simple groups.  But $H$ is normal in $N_2$ (since $\Phi_f(G_2)$ is normal in $G_2$). Thus $H$ is a direct product of non-abelian finite simple groups. 

This contradiction shows that $N\cap\Phi_f(G)=1$, which was the first assertion of the Lemma. The second assertion of the lemma now follows from Lemma \ref{frattini_under_epi}. This completes the proof. \qed\\[\baselineskip]

\section{Final comments}\label{sec7}

\subsection{An approach to Ivanov's question}

We will now discuss an approach to answering Ivanov's question 
(i.e. the nilpotency of $\Phi_f(G)$ for 
finitely generated subgroups $G$
of $\Gamma_g$) using the projective unitary representations described 
in \S 3. In the remainder of this section $G$ 
is an infinite, finitely generated subgroup of $\Gamma_g$.\\[\baselineskip]
The following conjecture is the starting point to this 
approach. Recall that a subgroup $G$ of $\Gamma_g$ is {\em reducible} if
there is a collection of essential simple closed curves $C$ on the
surface $\Sigma_g$, such that for any $\beta\in G$ there is a diffeomorphism
$\overline{\beta}:\Sigma_g\rightarrow \Sigma_g$ 
in the isotopy class of $\beta$ so that $\overline{\beta}(C)=C$. Otherwise
$G$ is called {\em irreducible}. As shown in \cite{Iv1} Theorem 2
an irreducible subgroup 
$G$ is either virtually an infinite cyclic subgroup
generated by a pseudo-Anosov element,
or, $G$ contains a free subgroup of rank $2$ generated
by a two pseudo-Anosov elements.\\[\baselineskip]
\noindent{\bf Conjecture:}~{\em If $G$ is a finitely generated irreducible 
non-virtually cyclic subgroup of $\Gamma_g$, then $\Phi_f(G)=1$.}\\[\baselineskip]
The motivation for this conjecture is that the irreducible (non-virtually 
cyclic)
hypothesis should be enough to guarantee that the image group
$\rho_p(\widetilde{G}) <\Delta_g$ is Zariski dense (with the same
adjoint trace-field).  Roughly speaking the irreducibility hypothesis
should ensure that there is no reason for Zariski density to fail (i.e. the
image is sufficiently complicated). Indeed, in this regard,
we note that an emerging theme in linear groups is that
random subgroups of linear groups are Zariski dense (see \cite{Ao} and
\cite{Ri} for example).  Below we discuss a possible approach to 
proving the Conjecture.

The idea now is to follow Ivanov's proof in \cite{Iv1} that the
Frattini subgroup is nilpotent. Very briefly if the subgroup 
is reducible then we first identify $\Phi_f$ on the pieces and then build
up to 
identify $\Phi_f(G)$. In Ivanov's 
argument, this involves passing to 
certain subgroups of $G$ (``pure subgroups''), understanding the Frattini
subgroup of these pure subgroups when restricted to the connected
components of $S\setminus C$, and then
building $\Phi(G)$ from this information. This uses several statements about the Frattini 
subgroup, at least one of which (Part (iv) of Lemma 10.2 of \cite{Iv1})
does not seem to easily extend to $\Phi_f$.\\[\baselineskip]
\noindent{\bf Remark:}~As a cautionary note to the previous discussion,
at present, it still remains conjectural that the image of a
fixed pseudo-Anosov element of $\Gamma_g$ under the representations
$\overline{\rho}_p$ is 
infinite order for big enough $p$
(which was raised in
\cite{AMU}).\\[\baselineskip]
\noindent{\bf An approach to the Conjecture:}\\[\baselineskip]
We begin by recalling that in \cite{Pl} Platonov
also proves that $\Phi_f(H)$ is nilpotent for
every
 finitely generated 
linear group $H$. Note that if $G$ is irreducible and virtually infinite cyclic
then 
$G$
is a linear group, and so \cite{Pl} implies that
$\Phi_f(G)$ is nilpotent.

Thus we now assume that $G$ is irreducible as in the conjecture.
Consider $\rho_p(\widetilde{\Phi_f(G)})$:
by Lemma \ref{frattini_under_epi} above
we deduce that $\rho_p(\widetilde{\Phi_f(G)})$ is a nilpotent 
normal subgroup of $\rho_p(\widetilde{G})$. Now
$\overline{\rho}_p(\Gamma_g)<\PSU(V_p,H_p;\BZ[\zeta_p])$ and it
follows from this that (in the notation of \S 3)
$\Delta_g < \Lambda_p=\SU(V_p,H_p;\BZ[\zeta_p])$. As
discussed in \cite{MR}, $\Lambda_p$ is a
cocompact arithmetic lattice in the algebraic group $\SU(V_p,H_p)$.
Thus $\rho_p(\widetilde{\Phi_f(G)})< \rho_p(\widetilde{G}) < \Lambda_p$. 
It follows from
general properties of cocompact lattices acting on symmetric spaces (see
e.g. \cite{Eb} Proposition 10.3.7) that $\rho_p(\widetilde{\Phi_f(G)})$ 
contains a
maximal normal abelian subgroup of finite index. Now there is a general bound
on the index of this abelian subgroup that is a function of the
dimension $N_g(p)$. However, in our setting, if the index can be bounded
by some fixed constant $R$ independent of $N_g(p)$, then 
we claim that $\Phi_f(G)$ can at least be shown to be 
finite. To see this we argue as follows.

Assume that $\Phi_f(G)$ is infinite.  Since $G$ is an irreducible subgroup
containing a free subgroup generated by a pair of pseudo-Anosov
elements, the same holds for the infinite normal subgroup $\Phi_f(G)$
(by standard dynamical properties of pseudo-Anosov elements, see for example
\cite{Lo} pp. 83--84).  

Thus we can find $x,y\in \Phi_f(G)$ a pair of non-commuting pseudo-Anosov
elements. Also note that $[x^t,y^t]\neq 1$ for all non-zero integers
$t$. From Lemma \ref{frattini_under_epi} we have that
$\rho_p(\widetilde{\Phi_f(G)}) < \Phi_f(\rho_p(\widetilde{G}))$ and
from the assumption above it therefore follows that
$\rho_p(\widetilde{\Phi_f(G)})$ contains a maximal normal abelian
subgroup $A_p$ of index bounded by $R$ (independent of $p$). Thus, setting
$R_1=R!$, we have
$[\overline{\rho}_p(x^{R_1}),\overline{\rho}_p(y^{R_1})]=1$ for all $p$. 
However, as noted above, $[x^{R_1},y^{R_1}]$ is a
non-trivial element of $G$, and by asymptotic faithfulness this cannot
be mapped
trivially for all $p$.  This is a contradiction.\\[\baselineskip]

\subsection{The profinite completion of $\Gamma_g$}

We remind the reader that the profinite completion $\widehat{\Gamma}$
of a group $\Gamma$ is the inverse limit of the finite quotients $\Gamma/N$ of
$\Gamma$. (The maps in the inverse system are the obvious ones: if $N_1
<N_2$ then $\Gamma/N_1\rightarrow \Gamma/N_2$.)
The Frattini subgroup $\Phi(G)$ of a profinite group $G$ is defined to be the
intersection of all maximal open subgroups of $G$.  Open subgroups
are of finite index, and if $G$ is finitely generated as a profinite group,
then Nikolov and Segal \cite{NS} show that finite index subgroups
are always open. Hence we can simply take $\Phi(G)$ to be
the intersection of all maximal finite index subgroups of $G$. 

Now if $\Gamma$ is a finitely generated residually
finite discrete group, the correspondence theorem between
finite index subgroups of $\Gamma$ and its profinite completion (see
\cite{RZ} Proposition 3.2.2) shows that 
$\overline{\Phi_f(\Gamma)} < \Phi(\widehat{\Gamma})$. 

There is a well-known connection between 
the center of a group $G$, denoted $Z(G)$ (profinite or otherwise), and
$\Phi(G)$. We include a proof for completeness. Note that for a profinite
group $\Phi(G)$ is a closed subgroup of $G$, $Z(G)$ is a closed subgroup
and by \cite{NS} the commutator subgroup $[G,G]$ is a closed subgroup.

\begin{lemma}
\label{profinitecenter}
Let $G$ be a finitely generated profinite group. Then 
$\Phi(G) > Z(G) \cap [G,G]$.\end{lemma}

\noindent{\bf Proof:}~Let $U$ be a maximal finite index subgroup of $G$, and
assume that $Z(G)$ is not contained in $U$. Then $<Z(G),U>=G$ by maximality.
It also easily follows that $U$ is a normal subgroup of $G$.   But then
$G/U = Z(G)U/U \cong Z(G)/(U\cap Z(G))$ which is abelian, and so $[G,G]<U$.
This being true for every maximal finite index subgroup $U$ we deduce that
$\Phi(G) >  Z(G) \cap [G,G]$ as required.\qed\\[\baselineskip]
We now turn to the following questions which were also 
part of the motivation
of this note. \\[\baselineskip]
\noindent{\bf Question 1:}~{\em For $g\geq 3$, is $Z(\widehat{\Gamma}_g)=1$?}

\medskip

\noindent{\bf Question 2:}~{\em For $g\geq 3$, is 
$Z(\widehat{\cal I}_g)=1$?}\\[\baselineskip]
Regarding Question 1, it is shown in \cite{HM} that 
the completion of $\Gamma_g$ arising from the congruence topology on $\Gamma_g$
has trivial center.  Regarding Question 2, if $Z(\widehat{\cal I}_g)=1$,
then the profinite topology on $\Gamma_g$ will induce
the full profinite topology on ${\cal I}_g$ (see \cite{LS} Lemma 2.6). 
Motivated by this and 
Lemma \ref{profinitecenter} we can also ask:\\[\baselineskip]
\noindent{\bf Question 1':}~{\em For $g\geq 3$, 
is $\Phi(\widehat{\Gamma}_g)=1$?}

\medskip

\noindent{\bf Question 2':}~{\em For $g\geq 3$, is 
$\Phi(\widehat{\cal I}_g)=1$?}\\[\baselineskip]

\noindent Although the results in this paper do not impact directly 
on Questions 1, 1',2
and 2', we note that since $\Gamma_g$ is finitely generated and perfect
for $g\geq 3$, it follows that $\widehat{\Gamma}_g$ is also perfect and hence 
$$
Z(\widehat{\Gamma}_g)<\Phi(\widehat{\Gamma}_g)$$ by Lemma \ref{profinitecenter}. As remarked above, the correspondence theorem gives 
$$\overline{\Phi_f(\Gamma_g)} < \Phi(\widehat{\Gamma}_g)~.$$ Thus our result that $\Phi_f(\Gamma_g)=1$ for $g\geq 3$ (which implies $\overline{\Phi_f(\Gamma_g)}=1$) is consistent with triviality of
$Z(\widehat{\Gamma}_g)$ (and similarly for $Z(\widehat{\cal I}_g)$).

\bigskip

\noindent Institut de Math{\'e}matiques de Jussieu-PRG (UMR 7586 du CNRS),\\
Equipe Topologie et G{\'e}om{\'e}trie Alg{\'e}briques,\\
Case 247, 4 pl. Jussieu,\\
75252 Paris Cedex 5, France.\\
\noindent Email:~gregor.masbaum@imj-prg.fr\\[\baselineskip] 
Department of Mathematics,\\
University of Texas\\
Austin, TX 78712, USA.\\
\noindent Email:~areid@math.utexas.edu

\end{document}